\documentclass{amsart}
\newcommand{\Z}{\mathbb{Z}}
\newcommand{\E}{\mathbf{E\,}}
\newcommand{\Var}{\mathbf{Var}\,}%
 \newtheorem{thm}{Theorem}[section]

 \theoremstyle{definition}
 
 \theoremstyle{remark}

 \numberwithin{equation}{section}

\begin{document}
%
%
%
%
%
%
%
%
%
\title{Sets with more differences than sums}
\author[J.-C. Schlage-Puchta]{Jan-Christoph Schlage-Puchta}

\address{
Building S22\\
Krijgslaan 281\\
9000 Gent\\
Belgium}

\email{jcsp@cage.ugent.be}



\date{February 24, 2009}
\begin{abstract}
We show that a random set of integers with density 0 has almost always
more differences than sums.
\end{abstract}

\maketitle

For a set $A\subseteq\Z$ set $A+A=\{a_1+a_2: A_i\in A\}$, and
$A-A=\{a_1-a_2: a_i\in A\}$. A finite set $A$ is called difference dominant,
if $|A-A|>|A+A|$, and sum dominant, if
$|A-A|<|A+A|$. Nathanson\cite{Nat} constructed infinite sequences of
sum dominant sets, and stated the opinion that the majority of all
subsets of $[1, n]$ is difference dominant. However, Martin and
O'Bryant\cite{MB} showed that the proportion of sum dominant sets is at
least $2\cdot 10^{-7}$. They conjectured that sets of density 0 are
almost always difference dominant. In this note we prove this
conjecture. More precisely, we have the following.

\begin{thm}
Let $p_n$ be a sequence of real numbers with $p_n\in[0, 1]$,
$p_n\rightarrow 0$ and $np_n\rightarrow\infty$. Let $\xi_{in}$, $0\leq
i\leq n-1$ be independent random variables satisfying
$P(\xi_{in}=1)=p_n$, and set $A_n=\{i:\xi_{in}=1\}$. Then the
probability that $A_n$ is difference dominant tends to 1.
\end{thm}

Martin and O'Bryant noted that for $p_n=o(n^{-3/4})$, this theorem
follows from the fact that in this case almost every set is a Sidon
set and has therefore almost twice as many differences as sums. 

\begin{proof}
We shall suppress the subscript $n$ throughout our argument. 

To simplify the computations we first deal with the case
$p=o(n^{-1/2})$. The number of elements of $A$ is asymptotically normal
distributed with mean and variance $np$, while
the expected number of solutions of the equation $x+y=u+v$ with $x, y,
u, v\in A$ is $\mathcal{O}(n^3p^4)$. Hence, with probability tending
to 1, we have
\begin{multline*}
|A-A|\geq |A|(|A|-1) - |\{(x, y, u, v\in A^4: x-y=u-v\}\\ >
(1-\epsilon)(np)^2 - \epsilon^{-1}n^3p^4
 > (1-2\epsilon)(np)^2 \geq
\frac{(|A|+1)|A|}{2}\geq |A+A|,
\end{multline*}
and our claim follows. Hence, from now on we shall assume that
$p>cn^{-1/2}$. 
Define random variables $\zeta_{1i}$, $\zeta_{2i}$ as
\[
\zeta_{1i} = \begin{cases} 1, & \exists a, b: a+b=i, \xi_a=\xi_b=1\\
0, & \mbox{otherwise}
\end{cases},\qquad
\zeta_{2i} = \begin{cases} 1, & \exists a, b: a-b=i, \xi_a=\xi_b=1\\
0, & \mbox{otherwise}
\end{cases},
\]
and set $S_j=\sum_{i\in\Z}\zeta_{ji}$. Then the probability of $A$ to be
difference dominant equals the probability of the event $S_2>S_1$. We
first compute the expectation 
of $S_j$. Noting that an even integer can be represented as the
sum of two different integers or as the double of an integer, we
obtain 
\begin{eqnarray*}
\E S_1 & = & \sum_{k=0}^{2n-2} \zeta_{1i}\\
 & = & (1-(1-p^2)^{\lfloor n/2\rfloor})) +
 2\underset{2\not| i}{\sum_{0\leq i\leq n-2}}(1-(1-p^2)^{(i-1)/2)}\\
 && \qquad + 2\underset{2| i}{\sum_{0\leq i\leq n-2}}(1-(1-p)(1-p^2)^{i/2-1})\\ 
 & = & 2\int_0^n 1-(1-p^2)^{t/2+\mathcal{O}(1)}\;dt+\mathcal{O}(1)\\
 & = & 2n-\frac{1-(1-p^2)^{n/2}}{-\log (1-p^2)^{1/2}} + \mathcal{O}(1)\\
 & = & 2n-\frac{2-2(1-p^2)^{n/2}}{p^2}+\mathcal{O}(1),
\end{eqnarray*}
and similarly
\begin{eqnarray*}
\E S_2 & = & \sum_{i=-n+1}^{n-1} \zeta_{2i}\\
 & = & 1 + 2\sum_{i=1}^{n-1}(1-(1-p^2)^{n-i}))+\mathcal{O}(1)\\
 & = & 2\int_0^n 1-(1-p^2)^t\;dt+\mathcal{O}(1)\\
 & = & 2n-\frac{1-(1-p^2)^n}{p^2}+\mathcal{O}(1).
\end{eqnarray*}
Since $p\gg n^{-1/2}$, we obtain $\E S_2 - \E S_1 \gg p^{-2}$.

Next, we give an upper bound for the variance of $S_j$. We have
\begin{eqnarray*}
\E S_j^2 & = & \big(\E S_j\big)^2 + 2\sum_{i<k}
\Big(P(\zeta_{ji}\zeta_{jk}=1)-P(\zeta_{ji}=1)P(\zeta_{jk}=1)\Big)\\
 &&\qquad + \sum_i P(\zeta_{ji}=1)-P(\zeta_{ji}=1)^2\\
 & = & \big(\E S_j\big)^2 +
\sum_i\Var\zeta_i + 2\sum_{i<k}
\Big(P(\zeta_{ji}\zeta_{jk}=1)-P(\zeta_{ji}=1)P(\zeta_{jk}=1)\Big)\\
\Var S_j & = & \sum_i\Var\zeta_i + 2\sum_{i<k}
\Big(P(\zeta_{ji}\zeta_{jk}=1)-P(\zeta_{ji}=1)P(\zeta_{jk}=1)\Big).\\
\end{eqnarray*}
Our aim is to show that $\Var S_j = o(p^{-4})$ for $j=1, 2$, our claim
then follows from Chebyshev's inequality together
with our estimate for $\E S_2-\E S_1$.

Obviously, the first term on the right-hand side is already of the
right magnitude, that is, it remains to bound the correlation of
$\zeta_{jk}$ and $\zeta_{ji}$ for $i<k$.

Clearly, the correlation of $\zeta_{ji}$ and $\zeta_{jk}$ is
non-negative, that is, it suffices to bound every summand from
above. We use two different estimates for
$P(\zeta_{ji}\zeta_{jk}=1)-P(\zeta_{ji}=1)P(\zeta_{jk}=1)$ depending
on whether $P(\zeta_{jk}=1)$ is close to 1 or not. First, we have
\[
P(\zeta_{ji}\zeta_{jk}=1)-P(\zeta_{ji}=1)P(\zeta_{jk}=1) \leq
P(\zeta_{ji}=1)\big(1-P(\zeta_{jk}=1)\big). 
\]
On the other hand, if $j=1$ and $i<k\leq n$, then
\begin{multline*}
P(\zeta_{1i}\zeta_{1k}=1)-P(\zeta_{1i}=1)P(\zeta_{1k}=1) \leq \\
P(\exists \mu, \nu, \kappa: \xi_{1\mu}=\xi_{1\nu}=\xi_{1\kappa}=1,
\mu+\nu=i, \mu+\kappa=k)\; \leq\; ip^3.
\end{multline*}
Similarly, if $i\leq n<k$, then
\[
P(\zeta_{1i}\zeta_{1k}=1)-P(\zeta_{1i}=1)P(\zeta_{1k}=1) \leq \max(0,
i+n-k) p^3.
\]
Hence, we have to show that the sum
\[
\sum_{i<k\leq n}\min(ip^3, \frac{(1-p^2)^k}{p^2})
\]
is of order $o(p^{-4})$. For each $k$ we either use the first or the
second estimate for all $i$, and obtain
\[
\sum_{k\leq n}\min\Big(k^2p^3, \frac{k(1-p^2)^k}{p^2}\Big) \ll \min_{k_0}
\Big(k_0^3p^3 + \frac{k_0(1-p^2)^{k_0}}{p^4}\Big).
\]
Putting $k_0=7p^{-2}\log p^{-1}$, the second term becomes $o(1)$,
while the first one is $\mathcal{O}(p^{-3}\log^4 p^{-1}) = o(p^{-4})$,
since $p\rightarrow 0$, which is of the desired size. A similar
computation shows that $S_2$ has variance $o(p^{-4})$, and we conclude
that the random variable $S_2-S_1$ has mean $p^{-2}$ and variance
$o(p^{-4})$, together with $p\rightarrow 0$ our claim follows.
\end{proof}


\begin{thebibliography}{9}
\bibitem{MB} G. Martin, K. O'Bryant, {\em Many sets have more sums than
  differences},  Additive combinatorics,  287--305, CRM Proc. Lecture
  Notes, 43, Amer. Math. Soc., Providence, RI, 2007.
\bibitem{Nat} M. B. Nathanson, {\em Sets with more sums than
  differences}, Integers {\bf 7} (2007), A5.
\end{thebibliography}
\end{document}